\documentclass{amsart}

\newtheorem{theorem}{Theorem}[section]
\newtheorem{lemma}[theorem]{Lemma}

\theoremstyle{definition}

\newtheorem{corollary}[theorem]{Corollary}
\newtheorem{conj}[theorem]{Conjecture}

\theoremstyle{remark}

\numberwithin{equation}{section}

\begin{document}

\title{The Ratio Monotonicity of the Boros-Moll Polynomials}

\author{William Y. C. Chen}
\address{Center for Combinatorics, LPMC-TJKLC
Nankai University
 Tianjin 300071, P. R. China}
\email{chen@nankai.edu.cn}
\thanks{ The authors are grateful to the referee
 for many helpful comments that lead to an improvement of an earlier
  version. We wish to thank
 Peter Paule for introducing us to
 this topic and for valuable discussions. This work was supported by  the 973
Project, the PCSIRT Project of the Ministry of Education, the
Ministry of Science and Technology, and the National Science
Foundation of China.}

\author{Ernest X. W. Xia}
\address{Center for Combinatorics, LPMC-TJKLC
Nankai University
 Tianjin 300071, P. R. China}
\email{xxw@cfc.nankai.edu.cn}


\date{June 26, 2008 and, in revised form, September 24, 2008.}



\begin{abstract}
In their study of a quartic integral, Boros and
Moll discovered a
special class of  Jacobi polynomials, which
we call the Boros-Moll
polynomials.
  Kauers and Paule proved the conjecture of
  Moll that these polynomials are
  log-concave. In this paper, we show that the Boros-Moll
  polynomials possess the ratio monotone
  property which implies the
  log-concavity and the spiral property. We conclude with a
  conjecture which is  stronger than Moll's conjecture on the
  $\infty$-log-concavity.
\end{abstract}

\maketitle

\section{Introduction}

In this paper, we aim to show that the Boros-Moll polynomials
satisfy the ratio monotone property which implies the log-concavity
and the spiral property.
 Boros and Moll
 \cite{George1999-1,George1999-2,George1999-3,George2001,George2004,Moll2002} explored  a special class of   Jacobi
polynomials in their study of a quartic integral. They have shown
that for any $a>-1$ and any nonnegative integer $m$,
\begin{equation}
 \int_{0}^\infty
\frac{1}{(x^4+2ax^2+1)^{m+1}}dx
=\frac{\pi}{2^{m+3/2}(a+1)^{m+1/2}}P_m(a),
\end{equation}
 where
\begin{align}
P_m(a)=\sum_{j,k}{2m+1 \choose 2j}{m-j \choose k}{2k+2j \choose
k+j}\frac{(a+1)^j(a-1)^k}{2^{3(k+j)}}.
\end{align}

Using Ramanujan's Master Theorem, Boros and Moll
\cite{George2001,Moll2002} derived the following formula
\begin{align}
P_m(a)=2^{-2m}\sum_{k}2^k{2m-2k \choose m-k} {m+k \choose k}(a+1)^k,
\end{align}
which indicates that the coefficients of $a^i$ in $P_m(a)$ is
positive for  $0\leq i \leq m$. Let $d_i(m)$ be  defined  by
\begin{equation} \label{pma}
P_m(a)=\sum_{i=0}^md_i(m)a^i.
\end{equation}
The polynomials $P_m(a)$ will be called the Boros-Moll polynomials,
and the sequence $\{d_i(m)\}_{0 \leq i \leq m}$ of the coefficients
will be called a Boros-Moll sequence. From \eqref{pma}, it follows
that
\begin{align}\label{Defi}
d_i(m)=2^{-2m}\sum_{k=i}^m2^k{2m-2k \choose m-k}{m+k \choose
k}{k\choose i}.
\end{align}
The readers can find in \cite{Am}  many proofs of this formula.
  Recall that $P_m(a)$ can be expressed as a
 hypergeometric function
\[
P_m(a)=2^{-2m}{2m \choose m}
{_2F_1}(-m,m+1;\frac{1}{2}-m;\frac{a+1}{2}),
\]
from which one sees that  $P_m(a)$ can be viewed as
 the    Jacobi
 polynomial  $P_m^{(\alpha, \beta)}(a)$ with
$\alpha=m+\frac{1}{2}$ and $\beta=-(m+\frac{1}{2})$, where
$P_m^{(\alpha,\beta)}(a)$ is given by
\[
P_m^{(\alpha,\beta)}(a)=\sum_{k=0}^m(-1)^{m-k}{m+\beta \choose m-k}
{m+k+\alpha+\beta \choose k}\left(\frac{1+a}{2}\right)^k.
\]

Boros and Moll \cite{George1999-2}  proved that the sequence
$\{d_i(m)\}_{0 \leq i \leq m}$ is unimodal and the maximum element
appears in the middle, namely,
\[d_0(m)<d_1(m)< \cdots
<d_{\left[\frac{m}{2}\right]}(m)>
d_{\left[\frac{m}{2}\right]+1}(m)>\cdots >d_m(m).\] They
 also
   established the unimodality by taking a
   different approach \cite{George1999-3}.
  Moll \cite{Moll2002}  conjectured that the sequence
  $\{d_i(m)\}_{0 \leq i \leq m}$ is
log-concave. Kauers and Paule \cite{Kauers2006} proved this
conjecture based on four recurrence relations found  using
 a
computer algebra approach. Two of these four recurrences have been
independently derived by Moll \cite{Moll2007} using the WZ-method.
Moreover, as will be seen, the two recurrences derived by Moll
easily imply the other two given by Kauers and Paule. These
recursions will be discussed in Section 2.

 Recall that a sequence $\{a_i\}_{0 \leq i \leq  m}$
 of positive numbers is said to be log-concave if
\[
\frac{a_0}{a_1} \leq \frac{a_1}{a_2} \leq \cdots  \leq
\frac{a_{m-1}}{a_m}.
\]
A polynomial is said to be log-concave if the sequence  of
 its coefficients is log-concave. It is easy to see that
if a sequence is log-concave then it is unimodal.
 A  sequence  $\{a_i\}_{0 \leq i \leq m}$
 of positive numbers is said to be spiral  if
\[ a_m \leq a_0 \leq a_{m-1} \leq
a_1 \leq \cdots \leq a_{\left[\frac{m}{2}\right]}.\]
 Similarly,  a
polynomial is said to be spiral if its sequence
 of coefficients is
spiral. It is easily seen that  a log-concave sequence is not
necessarily spiral, and vice versa. For example, $(2,10,3,1)$ is
spiral but not log-concave, whereas $(3,5,4,2,1)$ is lgo-concave but
not spiral. Chen and Xia \cite{Chen2007} discovered that the
$q$-derangement numbers are both spiral and log-concave, and
introduced the {\it ratio monotone} property defined below, which
implies both log-concavity and the spiral property. The
 purpose of this paper is to show that the Boros-Moll polynomials
possess the ratio monotone property.

  A  sequence $\{a_i\}_{0 \leq i \leq m}$ of
  positive numbers is said to be ratio monotone
 if
\begin{align}\label{D-1}
\frac{a_0}{a_{m-1}} \leq \frac{a_1}{a_{m-2}} \leq \cdots \leq
\frac{a_{i-1}}{a_{m-i}}\leq \frac{a_i}{a_{m-1-i}} \leq \cdots \leq
\frac{a_{\left[\frac{m}{2}\right]-1}}
{a_{m-\left[\frac{m}{2}\right]}} \leq  1
\end{align}
and
\begin{align}\label{D-2}
\frac{a_m}{a_0} \leq \frac{a_{m-1}}{a_{1}} \leq \cdots
 \leq \frac{a_{m-i}}{a_i} \leq  \frac{a_{m-1-i}}{a_{i+1}}
 \leq \cdots \leq \frac{a_{m-\left[\frac{m-1}{2}\right]}}
 {a_{\left[\frac{m-1}{2}\right]}}  \leq 1.
\end{align}
If every inequality relation   in \eqref{D-1} and \eqref{D-2}
becomes strict, we say that the sequence is strictly ratio monotone.
It is easy to see that the ratio monotonicity implies log-concavity.
In deeded, from \eqref{D-1} and \eqref{D-2}, we deduce that
\[
\frac{a_i}{a_{i-1}}\geq \frac{a_{m-1-i}}{a_{m-i}} \quad \mathrm{and}
\quad \frac{a_{i+1}}{a_i} \leq \frac{a_{m-1-i}}{a_{m-i}}.
\]
This gives
\[
\frac{a_i}{a_{i-1}} \geq \frac{a_{i+1}}{a_i}.
\]

The main result of this paper is stated as follows.

\begin{theorem}\label{Theo}
Let $m \geq 2$ be an integer. Then the Boros-Moll  sequence
$\{d_i(m)\}_{0 \leq i\leq m}$ satisfies the strictly ratio monotone
property. To be precise, we have
\begin{align}\label{Ratio-1}
\frac{d_m(m)}{d_0(m)}<\frac{d_{m-1}(m)}{d_1(m)}< \cdots
<\frac{d_{m-i}(m)}{d_i(m)}<\frac{d_{m-i-1}(m)}{d_{i+1}(m)}<\cdots
<\frac{d_{m-\left[\frac{m-1}{2}\right]}(m)}
{d_{\left[\frac{m-1}{2}\right]}(m)}<1
\end{align}
and
\begin{align}\label{Ratio-2}
\frac{d_0(m)}{d_{m-1}(m)}<\frac{d_1(m)}{d_{m-2}(m)}
<\cdots<\frac{d_{i-1}(m)}{d_{m-i}(m)}
<\frac{d_i(m)}{d_{m-i-1}(m)}<\cdots
<\frac{d_{\left[\frac{m}{2}\right]-1}(m)}
{d_{m-\left[\frac{m}{2}\right]}(m)}<1.
\end{align}
\end{theorem}

As a corollary of Theorem \ref{Theo}, we
  obtain  the spiral property of  the Boros-Moll sequences.
It is not clear whether there is a simpler way to
  verify this property directly.

\begin{corollary}
Let $m\geq 2$ be an  integer.
Then the Boros-Moll sequence
$\{d_i(m)\}_{0 \leq i \leq m}$ is spiral.
\end{corollary}

 The following example illustrates
 our main result.  For $m=8$, we have
 \allowdisplaybreaks
\begin{align*}
P_{8}(a)=&\frac{4023459}{32768}+\frac{3283533}{4096}a+
\frac{9804465}{4096}a^2+\frac{8625375}{2048}a^3
+\frac{9695565}{2048}a^4\\[6pt]
&+\frac{1772199}{512}a^5+\frac{819819}{512}a^6+\frac{109395}{256}a^7
+\frac{6435}{128}a^8.
\end{align*}
The strictly ratio monotone property is illustrated as follows:
\begin{align*}
 \frac {{\textstyle {\frac{ 6435}{ 128}}}}{\ \ \frac{ 4023459}{
32768}\ \  } &< \frac{\frac{109395}{256}} {\ \ \frac{3283533}{4096}\
\ }< \frac{ \frac{819819}{512}}{\ \ \frac{9804465}{4096}\ \  }<
\frac{\frac{1772199}{512}}{\ \ \frac{8625375}{2048}\ \ }<{\small 1},
\\[6pt]
\frac{\frac{4023459}{32768}}{\ \  \frac{109395}{256}\ \ }&< \frac{\
\  \frac{3283533}{4096}\  \ }{\frac{819819}{512}}< \frac{\ \
\frac{9804465}{4096}\ \  }{\frac{1772199}{512}}< \frac{\ \
\frac{8625375}{2048}\ \  }{\frac{9695565}{2048}}< 1.
\end{align*}
The spiral property of $P_8(x)$ is reflected by following order of
the coefficients:
\begin{align*}
\frac{6435}{128}<\frac{4023459}{32768}&<\frac{109395}{256}<
\frac{3283533}{4096}<\frac{819819}{512}\\[6pt]
&<\frac{9804465}{4096} <
\frac{1772199}{512}<\frac{8625375}{2048}<\frac{9695565}{2048}.
\end{align*}

 Based on the Moll conjecture on the $\infty$-log-concavity of the
 sequences $\{d_i(m)\}_{0\leq i\leq m}$, we conclude this paper with
 a stronger conjecture that these polynomials are infinitely ratio
 monotone. Numerical evidence seems to be supportive of this
 conjecture.

\section{Recurrence Relations}

We first  give a brief review of Kauers and Paule's approach
 to proving the log-concavity of  the Boros-Moll
sequence \cite{Kauers2006}. Our work employs the four recurrences
\begin{align}
d_i(m+1)=&\frac{m+i}{m+1}d_{i-1}(m)+\frac{(4m+2i+3)}{2(m+1)}d_i(m),
 \ \ \ \ 0 \leq i \leq m+1, \label{recu1}\\[6pt]
 d_{i}(m+1)=&\frac{(4m-2i+3)(m+i+1)}{2(m+1)(m+1-i)}d_i(m)
 \label{recu2} \\[6pt]
 & \quad -\frac{i(i+1)}{(m+1)(m+1-i)}d_{i+1}(m),
 \qquad \ \ \ 0 \leq i \leq
 m, \nonumber\\[6pt]
 d_i(m+2)=&\frac{-4i^2+8m^2+24m+19}{2(m+2-i)
 (m+2)}d_i(m+1)  \label{recu3} \\[6pt]
 & \quad -\frac{(m+i+1)(4m+3)(4m+5)}
 {4(m+2-i)(m+1)(m+2)}d_i(m), \qquad \ 0 \leq
  i \leq m+1,\nonumber
\end{align}
and for $0 \leq i \leq m+1$,
\begin{align}\label{recu4}
(m+2-i)(m+i-1)d_{i-2}(m)-(i-1)(2m+1)d_{i-1}(m)+i(i-1)d_i(m)=0.
\end{align}
These recurrences are derived  by Kauers and Paule \cite{Kauers2006}
with
 the RISC package MultiSum \cite{Wegsch1997}.
 In fact, the recurrences \eqref{recu3} and \eqref{recu4}
  are also derived independently by Moll \cite{Moll2007},
  and the other two relations (\ref{recu1}) and (\ref{recu2}) can
  be easily deduced
    from \eqref{recu3} and \eqref{recu4}.
    Based on the  four recurrence relations,
Kauers and Paule \cite{Kauers2006} used a computer algebra
 system  to
  derive the  next theorem, from which the log-concavity of the Boros-Moll
  sequence is  derived.

\begin{theorem} \label{tkp} For $0 < i< m$, we have
\begin{align}\label{condition}
d_i(m+1) \geq \frac{4m^2+7m+i+3}{2(m+1-i)(m+1)}d_i(m).
\end{align}
\end{theorem}

The   inequality  \eqref{condition} is also of vital importance for
our proof of the ratio monotonicity of the Boros-Moll sequences. We
note that the above inequality (\ref{condition}) is very tight. In
other words, the ratio
\[
\frac{(4m^2+7m+i+3)d_i(m)}{2(m+1-i)(m+1)d_i(m+1)},
\]
seems to be very close to $1$. For example, for $m=100$, the
smallest
 ratio is $0.998348$.

In order to establish the strict ratio monotonicity, we need a
slightly sharper version of \eqref{condition}. For example, we will
show that the inequality in (\ref{condition}) is strict for $1\leq
i\leq m-1$.

\begin{theorem}\label{improvement} Let $m \geq 2$. We have
\begin{align}\label{condition-1}
d_i(m+1)> \frac{4m^2+7m+i+3}{2(m+1-i)(m+1)}d_i(m), \qquad 1 \leq i
\leq m-1,
\end{align}
and
\begin{align}
d_0(m+1)&=\frac{4m+3}{2(m+1)}d_0(m), \label{I_1} \\[6pt]
d_m(m+1)&=\frac{(2m+3)(2m+1)}{2(m+1)}d_m(m) =
\frac{(2m+3)(2m+1)}{2(m+1)} 2^{-m}{2m\choose m}.  \label{I_2}
\end{align}
\end{theorem}

 To make this paper self-contained, we will present a detailed proof of the
 above improvement of Theorem \ref{tkp}.
 Before doing so, we remark that \eqref{recu3} and \eqref{recu4} can
  be also derived
 from \eqref{recu1} and \eqref{recu2}. Equating the right hand sides
  of \eqref{recu1} and \eqref{recu2} and replacing $i$ by $i-1$,
   we get \eqref{recu4}. Substituting  $i$ with $i+1$ and $m$ with
   $m+1$ in
  \eqref{recu1} and \eqref{recu2},
    respectively, we obtain
   two expressions for $d_{i+1}(m+1)$. This yields
    \allowdisplaybreaks
\begin{align}\label{F1}
 d_{i}(m+2)=&\frac{(4m-2i+7)(m+i+2)}{2(m+2)(m+2-i)}d_i(m+1)
  \\[6pt]
 & \quad -\frac{i(i+1)}{(m+2)(m+2-i)}
\left(\frac{m+i+1}{m+1}d_{i}(m)+\frac{(4m+2i+5)}{2(m+1)}
d_{i+1}(m)\right)\nonumber\\[6pt]
=&\frac{(4m-2i+7)(m+i+2)}{2(m+2)(m+2-i)}d_i(m+1)
-\frac{i(i+1)(m+i+1)}{(m+2)(m+2-i)(m+1)}d_i(m) \nonumber\\[6pt]
&\quad -\frac{i(i+1)(4m+2i+5)} {(m+2)(m+2-i)(2m+2)}
d_{i+1}(m).\nonumber
\end{align}
On the other hand, from \eqref{recu2}, we have
\begin{align}\label{F2}
d_{i+1}(m)=-\frac{(m+1)(m+1-i)}{i(i+1)}d_i(m+1)
+\frac{(m+i+1)(4m-2i+3)}{2i(i+1)}d_i(m).
\end{align}
Substituting \eqref{F2} into \eqref{F1}, we obtain \eqref{recu3}.

We now present a proof of Theorem \ref{improvement}.

\noindent {\it Proof.} Clearly, \eqref{I_1} follows from
\eqref{recu1} by setting  $i=0$, and  \eqref{I_2} can be obtained
from  \eqref{recu2} by setting  $i=m$.

We  proceed to prove \eqref{condition-1} by induction on $m$. It is
easy to verify that \eqref{condition-1} holds for $m=2$. We assume
that \eqref{condition-1} holds for $n \geq 2$, namely,
\begin{align}\label{assu}
d_i(n+1)> \frac{4n^2+7n+i+3}{2(n+1-i)(n+1)}d_i(n),
 \qquad 1 \leq i  \leq
n-1.
\end{align}
We aim to show that  \eqref{condition-1} holds for $n+1$, that is,
\begin{align} \label{n+1-in}
d_i(n+2) > \frac{4(n+1)^2+7(n+1)+i+3}{2(n+2)(n+2-i)}d_i(n+1),
 \ \ 1 \leq i \leq n.
\end{align}
Observe that for  $1 \leq i  \leq n-1$,
\begin{align*}
2(n+i&+1)(4n+3)(4n+5)(n+1-i)(n+1)
-2(4n^2+7n+i+3)\\[6pt]
&\times (n+1)(n+1-i)(4n+4i+5)=-4i(1+2i)(n+1)(n +1 -i)<0.
\end{align*}
 Hence we have  for $1 \leq i  \leq n-1$,
\begin{align}\label{F3}
\frac{4n^2+7n+i+3}{2(n+1-i)(n+1)}>\frac{(n+i+1)(4n+3)(4n+5)}
{2(n+1)(n+1-i)(4n+4i+5)}.
\end{align}
From the   inequalities  \eqref{F3} and \eqref{assu}, we  find that
for $1 \leq i \leq n-1$,
\begin{align}\label{assu-1}
d_i(n+1)>\frac{(n+i+1)(4n+3)(4n+5)} {2(n+1)(n+1-i)(4n+4i+5)}d_i(n).
\end{align}
It is easy to check that  {\small
\begin{align*}
 &\frac{{\displaystyle
\frac{(n+i+1)(4n+3)(4n+5)}{4(n+2-i)(n+1)(n+2)}}} { {\displaystyle
\frac{-4i^2+8n^2+24n+19}{2(n+2-i)(n+2)} -
\frac{4(n+1)^2+7(n+1)+i+3}{2(n+2-i)(n+2)}}}=\frac{(n+i+1)(4n+3)(4n+5)}
{2(n+1)(n+1-i)(4n+4i+5)}.
\end{align*}
} Hence  the inequality \eqref{assu-1} can be rewritten as
\begin{align*}
d_i(n+1)>\frac{\displaystyle{\frac{(n+i+1)(4n+3)(4n+5)}
{4(n+2-i)(n+1)(n+2)}}} {
{\displaystyle\frac{-4i^2+8n^2+24n+19}{2(n+2-i)(n+2)} -
\frac{4(n+1)^2+7(n+1)+i+3}{2(n+2-i)(n+2)}}}d_i(n).
\end{align*}
It follows that
\begin{align}\label{D-3}
&\frac{-4i^2+8n^2+24n+19}{2(n+2-i)(n+2)}d_i(n+1)
-\frac{(n+i+1)(4n+3)(4n+5)}
{4(n+2-i)(n+1)(n+2)}d_i(n) \\[6pt]
& \ \ \qquad \qquad
>\frac{4(n+1)^2+7(n+1)+i+3}
{2(n+2-i)(n+2)}d_i(n+1).\nonumber
\end{align}
From the recurrence  relation \eqref{recu3},  the left hand side of
\eqref{D-3} equals  $d_i(n+2)$. Thus we have verified  the
inequality (\ref{n+1-in}) for $1\leq i \leq n-1$. It is still
necessary to show that (\ref{n+1-in}) is true for $i=n$, that is,
\begin{equation}\label{i=n}
d_n(n+2)>\frac{4(n+1)^2+8n+10}{4(n+2)}d_n(n+1).
\end{equation}
Using the formula \eqref{Defi}, we get
 \begin{align*}
d_n(n+1)&=2^{-n-2}(2n+3){2n+2\choose n+1},\\[6pt]
d_n(n+2)&=\frac{(n+1)(4n^2+18n+21)}{2^{n+4}(2n+3)}{2n+4\choose n+2}.
 \end{align*}
It is easily checked that for $n \geq 1$,
 \begin{align*}
\frac{d_n(n+2)}{d_n(n+1)}=\frac{(n+1)(4n^2+18n+21)}{2(n+2)(2n+3)}
>\frac{4(n+1)^2+8n+10}{4(n+2)}.
\end{align*}
Hence the proof is complete by induction. \qed

\section{Preliminary Inequalities}

To prove the ratio  monotone property of the Boros-Moll polynomials,
we will establish the some inequalities based on the recurrence
relations derived by
 Kauers and
Paule \cite{Kauers2006} and Moll \cite{Moll2007}.

\begin{lemma}\label{Le}
Let $m \geq 2$  be an integer. Then we have
\begin{align}\label{ineqI}
\frac{m-j}{j+1} > \frac{d_{j+1}(m)}{d_j(m)}, \qquad
 1\leq j \leq
m-1.
\end{align}
\end{lemma}

\noindent {\it Proof.} From  \eqref{recu2} and Theorem
\ref{improvement},
  we  find  that  for $1 \leq j \leq m-1$,
  \begin{align*}
  (4m-2j+3)(m+j+1)d_{j}(m)-&2j(j+1)d_{j+1}(m)
  =2(m+1-j)(m+1)d_j(m+1)\\[6pt]
  &> (4m^2+7m+j+3)d_j(m),
  \end{align*}
 which implies \eqref{ineqI}. \qed

The following lemma gives an upper bound on the ratio ${d_i(m+1)/
d_i(m)}$, which is crucial for the proof of the main result of this
paper (Theorem \ref{Theo}).

\begin{lemma}\label{lemma1} Let $m \geq 2$ be a positive integer.
  We have  for $0 \leq i \leq m$,
 \begin{align}\label{in_1}
 d_i(m+1) \leq B(m,i)d_i(m),
 \end{align}
 where $B(m,i)$ is defined by
\begin{align}\label{denoteF}
B(m,i)=\frac{A(m,i)}{2(i+2)(4m+2i+5)(m+1)(m-i+1)}
\end{align}
with
\begin{align} \label{ANL}
A(m,i)=&\,30+96m^2+94m+37i+72m^2i+8m^2i^2-i^3
 \\[6pt]
&+99mi+5i^2+13mi^2+16m^3i+32m^3.\nonumber
\end{align}
\end{lemma}

\noindent {\it Proof.} We proceed by induction on $m$. It is easily
seen that the lemma  holds  for $m=2$. We assume that the lemma is
true for $n \geq 2$, i.e.,
\begin{align}\label{induction}
d_i(n+1) \leq B(n,i)d_i(n), \qquad 0 \leq i \leq n,
\end{align}
 where $B(n,i)$
is defined by \eqref{denoteF}. It will be shown that the lemma holds
  for $n+1$,   that is,
\begin{align}\label{ineq}
 d_i(n+2) \leq B(n+1,i)d_i(n+1),\ \ \ \  0
\leq i \leq n+1.
\end{align}
For $0\leq i \leq n$, let
\begin{align*}
F(n,i)=&(4n+2i+9)(i+2)(4n+5)(4n+3)(n+i+1),\\[6pt]
G(n,i)=&-2(-90-23i-202n+51i^3+60i^2-144n^2-32n^3\\[6pt]
&\quad -80n^2i-8n^2i^2-97ni+13ni^2-16n^3i+16ni^3+8i^4)(n+1).
\end{align*}
We claim  that
\begin{align} \label{inequality-1}
\frac{F(n,i)}{G(n,i)}\geq B(n,i), \qquad 0 \leq i \leq n.
\end{align}
Keeping in mind that $A(n,i)$ is defined by \eqref{ANL}, it is easy
to check that
\begin{align*}
2(i+2)&(4n+2i+5)(n+1)(n-i+1)F(n,i)-A(n,i)G(n,i)\\[6pt]
=&(128n^4i^4-32n^3i^5-80n^2i^6-16ni^7)
+(618n^3i^4-222ni^6-16i^7-284n^2i^5)\\[6pt]
&+(844ni^3-170i^4)+(1502n^2i^3-338i^5)+(984n^2i^4-142i^6)\\[6pt]
&+
(844n^3i^3-590ni^5)+256n^5i^2+720i+10i^3+788i^2+3984n^2i\\[6pt]
&+2656ni+3568ni^2+3136n^3i+4600n^3i^2+256n^5i\\[6pt]
&+1344n^4i+324ni^4+176n^4i^3+5908n^2i^2+1728n^4i^2.
\end{align*}
We are now in a position to see that the above expression is always
nonnegative since the expression in every parenthesis  is
nonnegative for $0 \leq i \leq n$. For example,
\[ 128n^4i^4-32n^3i^5-80n^2i^6-16ni^7 \geq
128n^4i^4-32n^4i^4-80n^4i^4-16n^4i^4 = 0.\] Thus we have
\begin{align}\label{Formula3}
2(i+2)(4n+2i+5)(n+1)(n-i+1)F(n,i)-A(n,i)G(n,i) \geq 0.
\end{align}
It is easy to see that $G(n,i) $ is positive for $0 \leq i\leq n$,
 and hence \eqref{inequality-1} can be deduced from \eqref{Formula3}.
 From the inductive hypothesis \eqref{induction} and
\eqref{Formula3}, it follows that
  for $0 \leq i \leq n$,
\begin{align}\label{T2}
\frac{F(n,i)}{G(n,i)}d_i(n)\geq B(n,i)d_i(n)\geq d_i(n+1).
\end{align}
It is a routine to verify that
\begin{align*}
\frac{(n+1+i)(4n+3)(4n+5)} {4(n+1)(n+2)(n+2-i)
\left({\displaystyle\frac{-4i^2+8n^2+24n+19}
{2(n+2-i)(n+2)}}-B(n+1,i)\right)}=\frac{F(n,i)}{G(n,i)}.
\end{align*}
From the above identity and \eqref{T2},  it follows that  for $0
\leq i \leq n$,
\begin{align}\label{T1}
&\frac{(n+1+i)(4n+3)(4n+5) d_i(n)} {4(n+1)(n+2)(n+2-i)
\left({\displaystyle \frac{-4i^2+8n^2+24n+19}
{2(n+2-i)(n+2)}}-B(n+1,i)\right)}\\[6pt]
&\qquad\qquad=\frac{F(n,i)}{G(n,i)}d_i(n)
 \geq d_i(n+1).\nonumber
\end{align}
Since
\[ \frac{-4i^2+8n^2+24n+19} {2(n+2-i)(n+2)}-B(n+1,i)\]
 is
positive for $0 \leq i \leq n$,  \eqref{T1}
 can be rewritten as
\begin{align}\label{fn}
&\frac{-4i^2+8n^2+24n+19} {2(n+2-i)
(n+2)}d_i(n+1) \\[6pt]
& \qquad \qquad
-\frac{(n+1+i)(4n+3)(4n+5)}{4(n+1)(n+2)(n+2-i)}d_i(n) \leq
B(n+1,i)d_i(n+1). \nonumber
\end{align}
From the recurrence relation \eqref{recu3}, we see that
\begin{equation}\label{fn-1}
\frac{-4i^2+8n^2+24n+19} {2(n+2-i)(n+2)}d_i(n+1)
-\frac{(n+1+i)(4n+3)(4n+5)}{4(n+1)(n+2)(n+2-i)}d_i(n)=d_i(n+2).
\end{equation}
In view of (\ref{fn}) and (\ref{fn-1}), we find that   the
inequality \eqref{ineq} holds for $0 \leq i \leq n$. It remains to
verify that (\ref{ineq}) holds for $i=n+1$, that is,
\begin{equation}\label{dn12}
 d_{n+1}(n+2) \leq B(n+1,n+1)d_{n+1}(n+1).
\end{equation}
By the definition (\ref{denoteF}) of $B(n,i)$, we have
\begin{align*}
B(n+1,n+1)=\frac{501+212n^3+692n^2+975n+24n^4}{2(n+3)(6n+11)(n+2)}.
\end{align*}
From the formula \eqref{Defi} for $d_i(m)$, we get
\[
d_{n+1}(n+1)=2^{-n-1}{2n+2 \choose n+1}
\]
and
\begin{align*}
d_{n+1}(n+2)=2^{-n-2}{2n+3 \choose n+1}+2^{-n-2}(n+2){2n+4 \choose
n+2}.
\end{align*}
Therefore, for $n \geq 0$, we have
\begin{align*}
\frac{d_{n+1}(n+2)}{d_{n+1}(n+1)}=\frac{(2n+3)(2n+5)}{2(n+2)} \leq
\frac{501+212n^3+692n^2+975n+24n^4}{2(n+3)(6n+11)(n+2)}.
\end{align*}
This completes the proof of the lemma. \qed

\begin{lemma}\label{Lem}
 Let $B(m,j)$ be  defined by  \eqref{denoteF} and $m\geq 2$ be an integer. Then we
have for $1 \leq j \leq m$,
\begin{align}\label{identity}
d_{j-1}(m)\leq \frac{2(m+1)B(m,j)-(4m+2j+3)}{2(m+j)}d_j(m).
\end{align}
\end{lemma}

\noindent {\it Proof.}  From the recurrence relation \eqref{recu1}
and Lemma \ref{lemma1}, we find that  for $0 \leq j \leq m$,
\begin{align}\label{identity3}
2(m+1)d_j(m+1)&=2(m+j)d_{j-1}
(m)+(4m+2j+3)d_j(m) \\[6pt]
 &\leq 2(m+1)B(m,j)d_j(m),\nonumber
\end{align}
where $B(m,j)$ is defined by  \eqref{denoteF}. Then
\eqref{identity3} implies   \eqref{identity}.
 \qed

\begin{lemma}\label{lemma2}
Let $m$ be a positive integer. For $ 0 \leq i \leq \frac{m}{2}$, we
have
\begin{align}\label{In_2}
\frac{2(2m-i)}{2(m+1)B(m,m-i)-(6m-2i+3)}>
\frac{2(m+1)B(m,i)-(4m+2i+3)}{2(m+i)},
\end{align}
where $B(m,i)$ is defined by \eqref{denoteF}.
\end{lemma}

\noindent {\it Proof.} For $0 \leq i \leq m$, let
\begin{align}
N(m,i)=\; &2(2m-i)(m-i+2)(6m-2i+5)(i+1), \label{NML}\\[6pt]
M(m,i)=\; &4(3m-i)(2m-i)(m-i)^2+(80m^3-155m^2i)
\label{MML}\\[6pt]
&\quad +\,(80m^2-108mi)+
(20m-20i)+(94mi^2-19i^3)+28i^2,\nonumber \\[6pt]
C(m,i)=\; &i(24m^2+52m+8m^2i+37mi+4i^3+12mi^2+20+19i^2+28i),
 \label{CML}\\[6pt]
D(m,i)=\; &2(i+2)(4m+2i+5)(m-i+1)(i+m). \label{DML}
\end{align}
Note that  $N(m,i), M(m,i), C(m,i)$ and $D(m,i)$ are all nonnegative
for $0 \leq i \leq \frac{m}{2}$, since  the sum in every parenthesis
in
 \eqref{NML}, \eqref{MML}, \eqref{CML} and \eqref{DML}
 is nonnegative for $0 \leq i \leq \frac{m}{2}$.
  It is easy
to check that
\begin{align*}
&N(m,i)D(m,i)-C(m,i)M(m,i)\\[6pt]
&\quad =(312m^5i^2+36m^2i^5+276m^3i^4-612m^4i^3-12mi^6)
+(2040m^4i^2-2533m^3i^3)
\\[6pt]
 &\qquad +(129mi^5-43i^6)+(384m^6-752m^5i)+(3568m^4-3328m^3i)\\[6pt]
 &\qquad +(1952m^5-2792m^4i)+(4280m^3i^2-2976m^2i^3)
 +(2800m^3-1240m^2i)\\[6pt]
 &\qquad +(3868m^2i^2-1080mi^3)+1240mi^2+1488mi^4+540i^4+800m^2
 +1159m^2i^4.
 \end{align*}
Observe that   the expression in every parenthesis in the above sum
is nonnegative for $0 \leq i \leq \frac{m}{2}$. Moreover, one sees
the term $800m^2$ is certainly positive. It follows that
 \begin{align}\label{formula5}
 N(m,i)D(m,i)-C(m,i)M(m,i)>0, \qquad 0 \leq i \leq \frac{m}{2}.
 \end{align}
Recall that $B(n,i)$ is defined by \eqref{denoteF}. It is easy to
check that
\begin{align*}
\frac{2(m+1)B(m,i)-(4m+2i+3)}{2(m+i)}&=\frac{C(m,i)}{D(m,i)},\\[6pt]
\frac{2(2m-i)}{2(m+1)B(m,m-i)-(6m-2i+3)}&=\frac{N(m,i)}{M(m,i)}.
\end{align*}
Thus the inequality \eqref{formula5} is equivalent to \eqref{In_2}.
This completes the proof of the lemma. \qed

\section{Proof of the Main Theorem}

Using the preliminary inequalities presented in the previous
section, we are ready to give a proof of Theorem \ref{Theo}.

\noindent {\it Proof.} It is clear that Theorem \ref{Theo}
  holds  for  $m=2,3,4$. We now assume that
$m \geq 5$. First we  consider \eqref{Ratio-1}. In order to verify
\begin{align}\label{I_3}
\frac{d_m(m)}{d_0(m)}< \frac{d_{m-1}(m)}{d_1(m)},
\end{align}
we invoke the formula \eqref{Defi} to get
\begin{align}\label{Formula1}
\frac{d_1(m)}{d_0(m)}=\frac{2^{-2m}\sum\limits_{k=1}^m2^k{2m-2k
\choose m-k}{m+k\choose m}k}{2^{-2m}\sum\limits_{k=0}^m2^k{2m-2k
\choose m-k}{m+k\choose m}} < \frac{\sum\limits_{k=1}^m2^k{2m-2k
\choose m-k}{m+k\choose m}m}{\sum\limits_{k=1}^m2^k{2m-2k \choose
m-k}{m+k\choose m}}=m,
\end{align}
and
\begin{align}\label{Formula2}
\frac{d_{m-1}(m)}{d_m(m)}=\frac{2^{-m}{2m-1 \choose m}+2^{-m}{2m
\choose m}m}{2^{-m}{2m \choose m}}>m.
\end{align}
 Combining \eqref{Formula1} and \eqref{Formula2}, we obtain
\[
 \frac{d_1(m)}{d_0(m)}< \frac{d_{m-1}(m)}{d_m(m)},
 \]
 which
 yields  \eqref{I_3}.

The next step is to  show that
\begin{align}\label{FA}
\frac{d_{m-i}(m)}{d_i(m)}<\frac{d_{m-i-1}(m)}{d_{i+1}(m)},
 \ \ \qquad 1 \leq i \leq \left[\frac{m-1}{2}\right]-1.
\end{align}
 By the assumption  $m \geq 5$,
 we have
$\left[\frac{m-1}{2}\right]-1 \geq 1$.
 Substituting $j$ with $i$ in \eqref{ineqI}, we have for
 $ 1\leq i
\leq \left[\frac{m-1}{2}\right]-1$,
\begin{align}\label{ineqI-1}
\frac{d_{i+1}(m)}{d_i(m)} < \frac{m-i}{i+1}.
\end{align}
On the other hand, since $ 1\leq i \leq
\left[\frac{m-1}{2}\right]-1$, we have $m-
\left[\frac{m-1}{2}\right] \leq m-i-1 \leq m-2$. Hence we may
substitute $j$ with $m-i-1$ in \eqref{ineqI} to deduce that
\begin{align}\label{ineqI-2}
\frac{d_{m-i-1}(m)}{d_{m-i}(m)} > \frac{m-i}{i+1}.
\end{align}
From \eqref{ineqI-1} and \eqref{ineqI-2}, it follows that for $
1\leq i\leq \left[\frac{m-1}{2}\right]-1$,
\[
\frac{d_{i+1}(m)}{d_i(m)} < \frac{m-i}{i+1}
<\frac{d_{m-i-1}(m)}{d_{m-i}(m)}.
\]
Hence we have verified \eqref{FA}.

It remains to show that  the last ratio
 in \eqref{Ratio-1} is smaller than $1$.
 Since $\left[\frac{m-1}{2}\right]<
m-\left[\frac{m-1}{2}\right]$, it is easily seen that for
$m-\left[\frac{m-1}{2}\right] \leq k \leq m$, we have
\[
{k\choose \left[\frac{m-1}{2}\right]} \geq {k\choose
m-\left[\frac{m-1}{2}\right]} .
\]
Based on the formula \eqref{Defi} and the above relation,
 we obtain
that
\begin{align*} d_{\left[\frac{m-1}{2}\right]}(m)
&=2^{-2m}\sum_{k=\left[\frac{m-1}{2}\right]}^m2^k{2m-2k\choose
m-k}{m+k\choose k}{k\choose \left[\frac{m-1}{2}\right]}\\[6pt]
&>2^{-2m}\sum_{k=m-\left[\frac{m-1}{2}\right]}^m2^k{2m-2k\choose
m-k}{m+k\choose k}{k\choose \left[\frac{m-1}{2}\right]}\\[6pt]
&\geq 2^{-2m}\sum_{k=m-\left[\frac{m-1}{2}\right]}^m2^k{2m-2k
\choose
m-k}{m+k\choose k}{k\choose m-\left[\frac{m-1}{2}\right]}\\[6pt]
&=d_{m-\left[\frac{m-1}{2}\right]}(m),
\end{align*}
leading to the relation
\[
\frac{d_{m-\left[\frac{m-1}{2}\right]}(m)}
{d_{\left[\frac{m-1}{2}\right]}(m)}<1.
\]
 This  completes the proof of \eqref{Ratio-1}.

We  now turn our attention to the proof  \eqref{Ratio-2}, which will
rely on the bound $B(n,i)$ and Lemmas \ref{Lem} and \ref{lemma2}.
First, rewrite \eqref{identity} as
\begin{align}\label{T_1}
\frac{d_{i-1}(m)}{d_{i}(m)}\leq
 \frac{2(m+1)B(m,i)-(4m+2i+3)}{2(m+i)},\quad 1\leq i \leq m.
\end{align}
For $1 \leq i\leq \left[\frac{m}{2}\right]$,  we have
$m-\left[\frac{m}{2}\right]\leq m-i\leq m-1$. It follows that
\begin{align}\label{Formula4}
2(m+1)&B(m,j)-(4m+2j+3)\\[6pt]
=&\frac{j(24m^2+8m^2j+52m+37mj+19j^2+28j+20+12mj^2+4j^3)}
{(j+2)(4m+2j+5)(m-j+1)},\nonumber
\end{align}
which is positive for $1 \leq j \leq m$. Substituting $j$ with $m-i$
in \eqref{Formula4}, we obtain that
\[
2(m+1)B(m,m-i)-(6m-2i+3)>0,  \qquad 1 \leq  i \leq
\left[\frac{m}{2}\right].
\]
Hence we can substitute $j$ with $m-i$ in  \eqref{identity} to
deduce  that for  $1 \leq i \leq \left[\frac{m}{2}\right]$,
\begin{align} \label{T_2}
\frac{d_{m-i}(m)}{d_{m-i-1}(m)}\geq
\frac{2(2m-i)}{2(m+1)B(m,m-i)-(6m-2i+3)}.
\end{align}
 Combining \eqref{T_1}, \eqref{T_2} and Lemma
\ref{lemma2}, we obtain that  for $1 \leq i \leq
\left[\frac{m}{2}\right]$,
\begin{align*}
 \frac{d_{i-1}(m)}{d_{i}(m)}  <
 \frac{d_{m-i}(m)}{d_{m-i-1}(m)},
\end{align*}
which can be restated as
\begin{align}\label{result}
\frac{d_{i-1}(m)}{d_{m-i}(m)}< \frac{d_i(m)}{d_{m-i-1}(m)}, \ \qquad
1 \leq i \leq \left[\frac{m}{2}\right].
\end{align}
At this point, it is necessary to  show that
 \begin{equation}\label{d-1} \frac{d_{\left[\frac{m}{2}\right]-1}(m)}
 {d_{m-\left[\frac{m}{2}\right]}(m)}<1.
 \end{equation}
 For $i=\left[\frac{m}{2}\right]$,  \eqref{result} becomes
  \begin{align}\label{identity4}
  \frac{d_{\left[\frac{m}{2}\right]-1}(m)}
 {d_{m-\left[\frac{m}{2}\right]}(m)} <
  \frac{d_{\left[\frac{m}{2}\right]}(m)}
 {d_{m-\left[\frac{m}{2}\right]-1}(m)}.
\end{align}
When $m$ is even, we have
$\left[\frac{m}{2}\right]=m-\left[\frac{m}{2}\right]$. From
\eqref{identity4}  it follows that
\[
\frac{d_{\left[\frac{m}{2}\right]-1}(m)}
 {d_{m-\left[\frac{m}{2}\right]}(m)}
 <\frac{d_{m-\left[\frac{m}{2}\right]}(m)}
 {d_{\left[\frac{m}{2}\right]-1}(m)},
\]
  which implies
    \eqref{d-1}. When $m$ is odd, we have
 $\left[\frac{m}{2}\right]=m-\left[\frac{m}{2}\right]-1$.
 Then \eqref{d-1} immediately follows from
 \eqref{identity4}.
   This  completes the proof
 of Theorem \ref{Theo}.   \qed

\section{A  Conjecture}

Moll made a conjecture on a property of the Boros-Moll sequences
which is stronger than the log-concavity. Given a sequence
$A=\{a_i\}_{0 \leq i \leq n}$,  define the  operator $\mathcal {L}$
 by $\mathcal{L}(A)=S=\{b_i\}_{0\leq i \leq n}$, where
 \[
b_i=a_i^2-a_{i-1}a_{i+1}, \qquad  0\leq i \leq n,
 \]
with the convention that $a_{-1}=a_{n+1}=0$.  We say that
$\{a_i\}_{0 \leq i \leq n}$ is $k$-log-concave if $\mathcal
{L}^j\left(\{a_i\}_{0 \leq i \leq n}\right)$ is log-concave for
every $0 \leq j \leq k-1$, and that $\{a_i\}_{0 \leq i \leq n}$ is
$\infty$-log-concave if $\mathcal {L}^k\left(\{a_i\}_{0 \leq i \leq
n}\right)$ is
 log-concave for every $k \geq 0$. Similarly, we say that
  $\{a_i\}_{0 \leq i \leq n}$ is
$j$-ratio-monotone
  (resp. $j$-strictly-ratio-monotone) if
  $\mathcal {L}^k\left(\{a_i\}_{0 \leq i \leq n}\right)$ is
 ratio monotone (resp. strictly ratio monotone)
 for every $0 \leq k \leq j-1$,
  and that
  $\{a_i\}_{0 \leq i \leq n}$ is
  $ \infty$-ratio-monotone
  (resp. $ \infty$-strictly-ratio-monotone)
  if $\mathcal {L}^k\left(\{a_i\}_{0 \leq i \leq n}\right)$ is
 ratio monotone (resp. strictly ratio monotone) for every $k\geq 0$.

 Moll \cite{Moll2002} has conjectured that the Boros-Moll
 sequence $\{d_i(m)\}_{0 \leq i \leq m}$
  is $\infty$-log-concave. We propose a stronger conjecture.

  \begin{conj}\label{conj}
  Suppose that $m\geq 2$ is a positive integer, then
   the Boros-Moll  sequence $\{d_i(m)\}_{0\leq i \leq m}$
  is
   $\infty$-strictly-ratio-monotone.
  \end{conj}

 We have verified that  the Boros-Moll sequence
 $\{d_i(m)\}_{0 \leq i \leq m}$ is $2$-strictly-ratio-monotone
  for $2 \leq m \leq 100$. For example,
  $\mathcal {L}\left(\{d_i(8)\}_{0 \leq i \leq 8}\right)$ is given by
  \begin{align*}
 b_0&=\frac{16188222324681}{1073741824},
\qquad b_1=\frac{46804848752277}{134217728},
\qquad b_2=\frac{39484127036475}{16777216},\\[6pt]
b_3&=\frac{53734360083525}{8388608}, \qquad
b_4=\frac{32860456870725}{4194304}, \qquad
b_5=\frac{4614148779669}{1048576},\\[6pt]
b_6&=\frac{284363773551}{262144}, \quad\qquad
 b_7=\frac{836466345}{8192},
\qquad\qquad \ \  b_8=\frac{41409225}{16384}.
  \end{align*}
It is easy to verify that
\begin{align*}
\frac{b_8}{b_0}<\frac{b_7}{b_1} <\frac{b_6}{b_2} <\frac{b_5}{b_3}
<1,\qquad  \frac{b_0}{b_7}<\frac{b_1}{b_6}
<\frac{b_2}{b_5}<\frac{b_3}{b_4}<1 .
\end{align*}

\bibliographystyle{amsplain}

\end{document}